\documentclass{amsart} 
\input{amssymb.sty} 
\title[Amenable actions and exactness for discrete groups] 
{Amenable actions and exactness \\ for discrete groups} 
\author{Narutaka Ozawa} 
\address{Department of Mathematical Science, 
University of Tokyo, Komaba, 153-8914, Japan } 
\email{narutaka@ms.u-tokyo.ac.jp} 
\thanks{The author is supported by JSPS} 
\date\today 
\newtheorem{thm}{Theorem} 
\newtheorem{lem}[thm]{Lemma} 
 
\theoremstyle{definition} 
\newtheorem{rem}[thm]{Remark} 
\newcommand{\C}{{\mathbb C}} 
\newcommand{\N}{{\mathbb N}} 
\newcommand{\M}{{\mathbb M}} 
\newcommand{\B}{{\mathbb B}} 
\newcommand{\V}{{\mathcal V}} 
\newcommand{\CP}{\mathrm{CP}} 
\newcommand{\hh}{{\mathcal H}} 
\newcommand{\id}{\mathrm{id}} 
\newcommand{\ran}{\mathop{\mathrm{ran}}} 
\newcommand{\lh}{\mathop{\mathrm{span}}} 

\begin{document} 
\begin{abstract} 
It is proved that a discrete group $G$ is exact if and only if 
its left translation action on the Stone-\v{C}ech compactification 
is amenable. 
Combining this with an unpublished result of Gromov, 
we have the existence of non exact discrete groups. 
\end{abstract} 
\maketitle  
In \cite{kw}, Kirchberg and Wassermann discussed exactness for groups. 
A discrete group $G$ is said to be exact if its reduced group 
$C^*$-algebra $C^*_\lambda(G)$ is exact. 
Throughout this paper, $G$ always means a discrete group 
and we identify $G$ with the corresponding convolution operators 
on $\ell_2(G)$. 

Amenability of a group action was discussed by 
Anantharaman-Delaroche and Renault in \cite{adr}. 
The left translation action of a group $G$ on its Stone-\v{C}ech 
compactification $\beta G$ was considered by Higson and Roe in \cite{hr}. 
This action is amenable if and only if the uniform Roe algebra 
$$UC^*(G):=C^*(\ell_\infty(G),G) 
=\overline{\lh}\{ s\ell_\infty(G) : s\in G\}\subset\B(\ell_2(G))$$ 
is nuclear. 
Since a $C^*$-subalgebra of an exact $C^*$-algebra is exact, 
$C^*_\lambda(G)$ is exact if $UC^*(G)$ is nuclear. 
In this article, we will prove the converse. 

A function $u\colon G\times G\to\C$ is called a positive definite kernel 
if the matrix $[u(s_i,s_j)]\in\M_n$ is positive 
for any $n$ and $s_1,\ldots,s_n\in G$. 
If $u$ is a positive definite kernel on $G\times G$ such that 
$u(s,s)\le1$ for all $s\in G$, then 
the Schur multiplier $\theta_u$ on $\B(\ell_2(G))$ defined by 
$$\theta_u(x)=[u(s,t)x_{s,t}]_{s,t\in G}$$ 
for $x=[x_{s,t}]\in \B(\ell_2(G))$ is 
a completely positive contraction. 
(See the section 3.6 in \cite{paulsen}.) 

\begin{lem}[Section 5 in \cite{paulsen}]\label{cp} 
Let $B$ be a $C^*$-algebra and $n\in\N$. 
Then, the map 
$$\CP(B,\M_n)\ni\phi\mapsto f_\phi\in(\M_n(B))^*_+$$ 
defined by 
$$f_\phi(X)=\sum_{i,j}\phi(x_{ij})_{ij}$$ 
for $X=[x_{ij}]\in\M_n(B)$ 
gives a bijective correspondence between 
the set of all completely positive maps 
$\CP(B,\M_n)$ from $B$ to $\M_n$ and the set of all positive linear 
functionals $(\M_n(B))^*_+$ on $\M_n(B)$. 
\end{lem} 
For vectors $\xi$ and $\eta$ in a Hilbert space $\hh$, 
we define a linear functional $\omega_{\xi,\eta}$ on $\B(\hh)$ by 
$\omega_{\xi,\eta}(x)=(x\xi,\eta)$ 
for $x\in\B(\hh)$. 
For a subset $\hh_0$ in $\hh$, we denote by $\V(\hh_0)$ 
the (possibly non-closed) linear span of 
$\{\omega_{\xi,\eta} : \xi,\eta\in\hh_0\}$. 

The following is a variation of Kirchberg's theorem. 
\begin{lem}\label{vec} 
Let $A$ be a unital exact $C^*$-subalgebra of $\B(\hh)$ and 
let $\hh_0$ be a total subset in $\hh$. 
Then, for any finite subset $E\subset A$ and $\varepsilon>0$, 
we have $\theta\colon A\to\B(\hh)$ such that 
\begin{enumerate} 
\item $\theta$ is of finite rank and unital completely positive, 
\item $\| \theta(x)-x \|<\varepsilon$ for all $x\in E$, 
\item there are $f_k\in\V(\hh_0)$ and operators $y_k$ in $\B(\hh)$ 
such that 
$$\theta(x)=\sum_{k=1}^df_k(x)y_k$$ 
for $x\in A$. 
\end{enumerate} 
\end{lem} 
\begin{proof} 
We may assume that $1\in E$. 
Since $A$ is exact, it follows from Kirchberg's theorem 
\cite{kirchberg, wassermann} that the inclusion map 
of $A$ into $\B(\hh)$ is nuclear. 
Thus, there are $n$ and unital completely positive maps 
$\phi\colon\B(\hh)\to\M_n$ and $\psi\colon\M_n\to\B(\hh)$ 
such that 
$$\|\psi\circ\phi(x)-x\|<\varepsilon/2\mbox{ for all }x\in E.$$ 
Let $f_\phi\in(\M_n(\B(\hh)))^*=\B(\hh^n)^*$ be 
the corresponding linear functional defined as in Lemma \ref{cp}. 
Since $\V(\hh_0^n)\cap\B(\hh^n)^*_+$ is weak$^*$ dense in $\B(\hh^n)^*_+$, 
we can approximate $f_\phi$ by linear functionals in 
$\V(\hh_0^n)\cap\B(\hh^n)^*_+$ in the weak$^*$-topology. 
It follows that we can approximate $\phi$ in the point-norm topology 
by completely positive maps $\phi'$ such that $\phi'(\cdot)_{ij}$ is 
in $\V(\hh_0)$. 
Thus, for arbitrary $0<\delta<1$, we may find such $\phi'$ with 
$$\|\phi'(x)-\phi(x)\|<\delta\mbox{ for all }x\in E.$$  
Let $p=\phi'(1)$. 
Since $1\in E$ and $0<\delta<1$, $p$ is invertible. 
Thus, we can define a unital completely positive map 
$\phi''\colon A\to\M_n$ by 
$$\phi''(x)=p^{-1/2}\phi'(x)p^{-1/2}$$ 
for $x\in A$. 
Taking $\delta>0$ sufficiently small, we have 
$$\|\phi''(x)-\phi(x)\|<\varepsilon/2\mbox{ for all }x\in E.$$ 
Finally, put $\theta=\psi\circ\phi''$ and we are done. 
\end{proof} 
The following was inspired by the work of Guentner and Kaminker \cite{gk}. 
\begin{thm}\label{main}  
For a discrete group $G$, the following are equivalent. 
\begin{enumerate} 
\item The reduced group $C^*$-algebra $C^*_\lambda(G)$ is exact. 
\item For any finite subset $E\subset G$ and any $\varepsilon>0$, 
there are a finite subset $F\subset G$ and a positive definite 
kernel $u\colon G\times G\to\C$ such that 
$$u(s,t)\neq0\mbox{ only if }st^{-1}\in F$$ 
and that 
$$|1-u(s,t)|<\varepsilon\mbox{ if }st^{-1}\in E.$$ 
\item The uniform Roe algebra $UC^*(G)$ is nuclear. 
\end{enumerate} 
\end{thm} 
\begin{proof} 
(i)$\Rightarrow$(ii). 
We follow the proof of Theorem 3.1 in \cite{gk}. 
We give ourselves a finite subset $E\subset G$ and $\varepsilon>0$. 
By the assumption, there is $\theta\colon C^*_\lambda(G)\to\B(\ell_2(G))$ 
satisfying the conditions in Lemma \ref{vec} for 
$E\subset\C^*_\lambda(G)$, $\varepsilon$ and $\hh_0=\{\delta_s : s\in G\}$. 
If 
$$\theta=\sum_{k=1}^d\omega_{\delta_{p(k)},\delta_{q(k)}}\otimes y_k$$ 
for $p(k),q(k)\in G$, then we put $F=\{q(k)p(k)^{-1} : k=1,\ldots,d\}$. 

We define $u\colon G\times G\to\C$ by 
$$u(s,t)=(\delta_s,\theta(st^{-1})\delta_t)$$ 
for $s,t\in\Gamma$. 
It is not hard to check that $u$ has the desired properties. 

(ii)$\Rightarrow$(iii). 
Let $\{ E_i\}_{i\in I}$ be an increasing net of finite subsets of $G$ 
containing the unit $e$. 
By the assumption, there are finite subsets $F_i\subset G$ and 
a net of positive definite kernels $\{ u_i\}_{i\in I}$ 
such that 
$$u_i(s,t)\neq0\mbox{ only if }st^{-1}\in F_i$$ 
and that 
$$|1-u_i(s,t)|<|E_i|^{-1}\mbox{ if }st^{-1}\in E_i.$$ 
We may assume that $u(s,s)\le 1$ for all $s\in G$. 
Let $\theta_i\colon\B(\ell_2(G))\to\B(\ell_2(G))$ be the Schur multiplier 
associated with the positive definite kernel $u_i$. 
Then, $\theta_i$'s are completely positive contractions and 
it can be seen that 
$$\ran\theta_i\subset\lh\{ s\ell_\infty(G) : s\in F_i\}\subset UC^*(G)$$ 
and that 
$$\lim_{i\in I}\|\theta_i(x)-x\|=0\mbox{ for all }x\in UC^*(G).$$ 
Let $\Phi$ be the restriction map from $\B(\ell_2(G))$ onto 
$\ell_\infty(G)$, i.e., $\Phi(x)(s)=(x\delta_s,\delta_s)$ 
for $x\in\B(\ell_2(G))$. 
For $i\in I$ and $s\in F_i$, we define a complete contraction 
$\sigma_i^s\colon\B(\ell_2(G))\to\ell_\infty(G)$ by 
$$\sigma_i^s(x)=\Phi(s^{-1}\theta_i(x))$$ 
for $x\in\B(\ell_2(G))$. 
Then, we have $\theta_i(x)=\sum_{s\in F_i}s\sigma_i^s(x)$ 
for $i\in I$ and $x\in\B(\ell_2(G))$. 

To prove that $UC^*(G)$ is nuclear, we take a unital $C^*$-algebra $B$. 
It suffices to show that $UC^*(G)\otimes_{\min}B=UC^*(G)\otimes_{\max}B$. 
Let 
$$Q\colon UC^*(G)\otimes_{\max}B\to UC^*(G)\otimes_{\min}B$$ 
be the canonical quotient map. 
Since $\ell_\infty(G)$ is nuclear, we observe that 
$$\sigma_i^s\otimes\id_B\colon UC^*(G)\otimes_{\min}B\to 
\ell_\infty(G)\otimes_{\min}B\subset UC^*(G)\otimes_{\max}B$$ 
is a well-defined contraction. 
Since $\theta_i$'s are completely positive contractions, so are 
$$\theta_i\otimes\id_B\colon 
UC^*(G)\otimes_{\max}B\to UC^*(G)\otimes_{\max}B$$ 
(see Theorem 10.8 in \cite{paulsen}) and we have 
$$\lim_{i\in I}\|\theta_i\otimes\id_B(z)-z\|=0\mbox{ for all }
z\in UC^*(G)\otimes_{\max}B.$$ 
Combining this with the factorization 
$$\theta_i\otimes\id_B(z)= 
\sum_{s\in F_i}(s\otimes 1)(\sigma_i^s\otimes\id_B(Q(z))),$$ 
we see that $\ker Q=\{ 0\}$. 

(iii)$\Rightarrow$(i). This is obvious. 
\end{proof} 
\begin{rem} 
Recently, Gromov found examples of finitely presented groups 
which are not uniformly embeddable into Hilbert spaces \cite{gromov}. 
As it was suggested in \cite{gk}, 
these examples of Gromov are indeed not exact 
since the condition (ii) in Theorem \ref{main} assures uniform embeddings. 

Our definition of exactness for discrete groups is different from 
the original one \cite{kw}, but this is justified by Theorem 5.2 in \cite{kw}. 
Also, we can reprove this using Theorem \ref{main}. 
Indeed, for  a closed 2-sided ideal $I$ in a $C^*$-algebra $A$, 
the corresponding sequence 
$$0\longrightarrow UC^*(G;I)\longrightarrow UC^*(G;A) 
\longrightarrow UC^*(G;A/I)\longrightarrow 0$$
is exact if the condition (ii) in Theorem \ref{main} holds, where 
for a $C^*$-algebra $B\subset\B(\hh)$, we set 
$$UC^*(G;B)=C^*(\ell_\infty(G;B), \C1_\hh\otimes\lambda(G)) 
\subset\B(\hh\otimes_2\ell_2(G)).$$ 

Our definition of amenable action is different from 
Definition 2.2 in \cite{hr} and there is a delicate problem 
when we are dealing with non second countable space, 
but this is justified when $G$ is countable. 
Indeed, let $u$ be a positive definite kernel 
as in the condition (ii) in Theorem \ref{main}. 
We may assume that $u(s,s)=1$ for all $s\in G$. 
Regarding $u$ as a positive element in $UC^*(G)$, 
we let $\xi_s=u^{1/2}\delta_s\in\ell_2(G)$. 
Then, we have $(\xi_t,\xi_s)=u(s,t)$ for $s,t\in G$. 
Now, define $\mu\colon G\to\ell_1(G)$ by 
$\mu_s(t)=|\xi_s(t^{-1}s)|^2$ for $s,t\in G$. 
It can be verified that $\|\mu_s\|_1=\|\xi_s\|_2^2=1$ and that 
$$\| s\mu_t-\mu_{st}\|_1=\|\,|\xi_t|^2-|\xi_{st}|^2\,\|_1 
\le\|\,|\xi_t|-|\xi_{st}|\,\|_2\,\|\,|\xi_t|+|\xi_{st}|\,\|_2\\ 
\le 2\sqrt{2\varepsilon}$$ 
for all $s\in E$ and $t\in G$. 
We observe that the range of $\mu$ is relatively weak$^*$ compact 
in $\ell_1(G)$ since $u^{1/2}$ is in $UC^*(G)$. 
Thus, we can extend $\mu$ on the Stone-\v{C}ech compactification $\beta G$ 
by continuity. 
This completes the proof of our claim. 

 From the theory of exact operator spaces developed by Pisier \cite{pisier}, 
it is enough to check the condition (ii) in Theorem \ref{main} for 
finite subsets $E$ that contain the unit $e$ and are contained in 
a given set of generators of $G$. 

See \cite{yu, higson, hr, gk} for the connection to the Novikov conjecture. 
\end{rem} 
\providecommand{\bysame}{\leavevmode\hbox to3em{\hrulefill}\thinspace} 

\end{document}